\documentclass[11pt]{amsart}
\usepackage{amsmath,amsfonts,amsthm,amsopn,amssymb}
\usepackage{cite,marginnote}
\usepackage{amsmath}
\usepackage{color,enumitem,graphicx}
\usepackage[colorlinks=true,urlcolor=blue,
citecolor=red,linkcolor=blue,linktocpage,pdfpagelabels,
bookmarksnumbered,bookmarksopen]{hyperref}
\usepackage[english]{babel}
\usepackage[left=2.8cm,right=2.8cm,top=3cm,bottom=3cm]{geometry}

\usepackage[hyperpageref]{backref}

\numberwithin{equation}{section}

\makeindex
\newtheorem{theorem}{Theorem}[section]
\newtheorem{proposition}[theorem]{Proposition}
\newtheorem{lemma}[theorem]{Lemma}
\newtheorem{remark}[theorem]{Remark}
\newtheorem{example}[theorem]{Example}
\newtheorem{corollary}[theorem]{Corollary}
\newtheorem{definition}[theorem]{Definition}

\newcommand{\RN}{\mathbb R^N}

\newcommand{\iy}{\infty}

\newcommand{\s}{\section}

\newcommand{\DD}{\Delta}
\newcommand{\g}{\gamma}
\newcommand{\G}{\Gamma}
\newcommand{\na}{\nabla}
\newcommand{\la}{\lambda}

\newcommand{\R}{\mathbb R}

\newcommand{\ti}{\tilde}

\newcommand{\rhu}{\rightharpoonup}
\newcommand{\lan}{\langle}
\newcommand{\ran}{\rangle}
\newcommand{\e}{\varepsilon}
\newcommand{\vp}{\varphi}

\newcommand{\lab}{\label}
\newcommand{\bt}{\begin{theorem}}
\newcommand{\et}{\end{theorem}}
\newcommand{\bl}{\begin{lemma}}
\newcommand{\el}{\end{lemma}}
\newcommand{\bd}{\begin{definition}}
\newcommand{\ed}{\end{definition}}
\newcommand{\bc}{\begin{corollary}}
\newcommand{\ec}{\end{corollary}}
\newcommand{\bp}{\begin{proof}}
\newcommand{\ep}{\end{proof}}
\newcommand{\bx}{\begin{example}}
\newcommand{\ex}{\end{example}}
\newcommand{\bi}{\begin{exercise}}
\newcommand{\ei}{\end{exercise}}
\newcommand{\bo}{\begin{proposition}}
\newcommand{\eo}{\end{proposition}}
\newcommand{\br}{\begin{remark}}
\newcommand{\er}{\end{remark}}
\newcommand{\be}{\begin{equation}}
\newcommand{\ee}{\end{equation}}
\newcommand{\ba}{\begin{align}}
\newcommand{\ea}{\end{align}}
\newcommand{\bn}{\begin{enumerate}}
\newcommand{\en}{\end{enumerate}}
\newcommand{\bg}{\begin{align*}}
\newcommand{\eg}{\end{align*}}
\newcommand{\bcs}{\begin{cases}}
\newcommand{\ecs}{\end{cases}}

\newcommand{\NN}{{\mathbb N}}

\newcommand{\bean}{\begin{eqnarray*}}
\newcommand{\eean}{\end{eqnarray*}}

\title[ground states]{Ground states of nonlinear fractional Schr\"{o}dinger equation involving critical growth}

\author[H.\ Jin]{Hua Jin}
\author[W. B.\ Liu]{Wenbin Liu}

\address[H.\ Jin]{\newline\indent College of Science
\newline\indent
China University of Mining and Technology
\newline\indent
Xuzhou, 221116, China}
\email{\href{mailto:huajin@cumt.edu.cn}{huajin@cumt.edu.cn}}

\address[W. B.\ Liu]{\newline\indent College of Science
\newline\indent
China University of Mining and Technology
\newline\indent
Xuzhou, 221116, China}
\email{\href{mailto:liuwenbin-xz@163.com}{liuwenbin-xz@163.com}}

\thanks{W. Liu is the corresponding author. Telephone number: (86-516) 83591530.}

\date{\today}

\begin{document}

\begin{abstract}
In this paper, we are concerned with the ground state solutions of nonlinear fractional Schr\"{o}dinger equation involving critical growth. Without Ambrosetti-Rabinowitz condition and monotonicity condition on the nonlinearity, we get the existence of ground state solutions of such equation when the potential $V(x)$ is not a constant and not radial.
\end{abstract}

\maketitle

\s{Introduction}
\renewcommand{\theequation}{1.\arabic{equation}}

\
The fractional Laplacian $(-\Delta)^s$ is a classical linear integro-differential operator of order $s$. The main feature, and also its main difficulty, is that it is a non-local operator. Recently, a great deal of attention has been devoted to the fractional Laplacian and non-local operators of elliptic type, both for their interesting theoretical structure and concrete applications. The fractional Laplacian $(-\Delta)^s$ arises in the description of various phenomena in the applied science, such as the thin obstacle problem \cite{L.A.Caffarelli,Silvestre}, phase transition \cite{G.Alberti,Y.Sire}, Markov processes \cite{Q.Y.Guan} and fractional quantum mechanics \cite{Laskin} and the references therein for more details.

The fractional Schr\"{o}dinger equation formulated by Nick Laskin \cite{NLaskin1,NLaskin2,Laskin} has the following form
\begin{align}
\label{fractionequation}
i \vp_t-(-\Delta)^s\varphi-V(x)\varphi+f(\varphi)=0,\ \ \ (x,t)\in \mathbb{R}^N\times\mathbb{R},
\end{align}
where $s\in (0,1)$, $N>2s$, $\vp$ is the  wavefunction and $V(x)$ is the potential energy. The fractional
quantum mechanics has been discovered as a result of expanding the Feynman path integral, from
the Brownian-like to the L\'{e}vy-like quantum mechanical paths. Since we are concerned with the standing wave solutions with the form
\begin{align*}
\vp(x,t)=e^{-iwt}u(x), \ \ w\in \mathbb{R},
\end{align*}
then equation (\ref{fractionequation}) can be converted into
\begin{align}
\label{p1}
(-\Delta)^su+V(x)u=f(u), \ \ x\in \RN.
\end{align}

When $s=1$, equation (\ref{p1}) gives back to the classical nonlinear Schr\"{o}dinger equation
\be
\lab{classicalequation}
-\Delta u+V(x)u=f(u), \ \ x\in \RN,
\ee
which has been studied theoretically and numerically in the last decades. We should emphasize that the potential $V(x)$ plays a crucial role concerning the existence of nontrivial solutions and the existence of ground state solutions.  If the potential $V(x)$ is a constant, namely (\ref{classicalequation}) is autonomous, in the celebrated paper \cite{Berestycki}, H. Berestycki and P. Lions first proposed the Berestycki-Lions conditions which are almost optimal for the existence of  ground state solutions in the subcritical case. The authors investigated the constraint minimization problem and use the Schwarz symmetrization in  $H_r^1(\RN)$. For the critical nonlinearity $f$, because of the lack of compactness of $H^1(\mathbb{R}^N)\hookrightarrow L^{2^*}(\mathbb{R}^N)$, the existence of ground state solutions of problem (\ref{classicalequation}) becomes rather more complicated. In \cite{zjjzwm2}, the critical case was considered by modifying the minimization methods with constrains. Since the radial symmetry plays a crucial role, the method is invalid for non-radial case.

In the  non-autonomous case, that is $V(x)\not\equiv V$, where $V$ is a constant, the main obstacle to get the existence of solutions or ground state solutions is the boundedness of the Palais-Smale (PS for short) sequence because of no some global conditions on $f$, such as the Ambrosetti-Rabinowitz (A-R for short) condition. Moreover, the lack of compactness due to the unboundedness of the domain prevents us from checking the PS condition.  To avoid the difficulties mentioned above, in the seminal paper \cite{L.JeanjeanandKazunagaTanaka}, Jeanjean and Tanaka used an indirect approach developed in \cite{L.Jeanjean} to get a bounded PS (BPS for short) sequence for the energy functional $I$, then the existence of  positive solutions and moreover ground state solutions is obtained in the subcritical case when nonlinearity $f$ and potential $V(x)$ satisfy the conditions as follows
\begin{itemize}

\item [$(F_1)$] $f\in C(\R^+,\R)$, $f(0)=0$ and $f'(0)$ defined as   $\lim_{t\to 0^+}f(t)/t$ exists,

\item [$({F_2})$]  there is $p<\iy$ if $N=2$, $p<2^*-1$ if $N\geq 3$ such that $\lim_{t\to \iy}{f(t)/t^{-p}}=0$,

\item [$({F_3})$]  $\lim_{t\to \iy}{f(t)/t}=+\iy$,

\item [$(V_1')$] $f'(0)<\inf\sigma(-\DD+V(x))$, where $\sigma(-\DD+V(x))$ denotes the spectrum of the self-adjoint operator $-\DD+V(x)$,
\item [$(V_2')$] $V\in C(\RN,\R)$, $V(x)\to V(\iy)\in \R$ as $|x|\to\iy$,
\item [$(V_3')$] $V(x)\leq V(\iy)$,
\item [$(V_4')$] there exists a function $\phi\in L^2(\RN)\cap W^{1,\iy}(\RN)$ such that
\begin{align*}
|x||\na V(x)|\leq \phi^2(x), \forall x\in\RN.
\end{align*}
\end{itemize}
Here, the decay condition $(V_4')$ is crucial to derive the boundedness of the PS sequence. For the critical case, the problem is  different and more difficulty. In \cite{J.Zhang}, by use of the indirect approach developed in \cite{L.Jeanjean}, the authors completed the proof of existence of ground state solutions in the critical case with the same conditions on $V(x)$. As for nonlinearity $f$, the following conditions are satisfied
\begin{itemize}

\item [$(F_1')$] $f\in C(\R^+,\R)$, $f(t)=o(t)$ as $t\to o^+$,

\item [$({F_2'})$]  $\lim_{t\to +\iy}f(t)/{t^{2^*-1}}=K>0$, where $2^*=\frac{2N}{N-2}$,

\item [$({F_3'})$]  there exist $D>0$ and $2<q<2^*$ such that $f(t)\geq Kt^{2*-1}+Dt^{q-1},\forall t\geq 0$,

\item [$({F_4'})$] $f\in C^1(\R^+,\R),|f'(t)|\leq C(1+|t|^{\frac{4}{N-2}})$.
\end{itemize}

Now, we return our attention to the fractional and non-local problems. With the aid of the extended techniques developed
by Caffarelli and Silvestre \cite{Caffarelli}, some existence and nonexistence of Dirichlet problems involving the fractional Laplacian on bounded domain have been established, see \cite{Barriosa1,Cabre,J.Tan} and so on. For the general fractional Schr\"{o}dinger equation
\begin{align*}
(-\Delta)^su+V(x)u=f(x,u),x\in\RN
\end{align*}
in subcritical or critical case,  many results have been obtained on the existence of ground state solutions, positive solutions, multiplicity of standing wave solutions, symmetry of solutions and so forth, under the different conditions on $V(x)$ and $f$, for example monotonicity condition, A-R condition, see \cite{B.Barriosa,bisci,Caffarelli,changxiaojun,chenzheng,tengkaiming,weiyuanhong,zhangxia} and the references therein.

As is well known, the existence and concentration phenomena of solutions on the singularly perturbed fractional Schr\"{o}dinger equation
\begin{align*}
\varepsilon^{2s}(-\Delta)^su+V(x)u=f(u), x\in \RN
\end{align*}
is also a hot topic. For this subject we refer, for example, to \cite{Alvesfrac,Davila,FMM,xiaominghe,seok,shangxzhangj} and the references therein.

Now, let us say more about the existence of  ground state solutions of a class of fractional scalar field equations
\begin{align*}
(-\Delta)^su+V(x)u=f(u),x\in\RN.
\end{align*}
When $f(u)-V(x)u=g(u)$, the authors \cite{changxiaojun2} obtained the existence of radial positive ground state solutions under the general Berestycki-Lions type assumptions in the case of subcritical growth. By using  the fractional Pohoz\v{a}ev identity and the monotonicity trick of Struwe-Jeanjean, they showed that the compactness still holds under their assumptions without the Strauss type radial lemma in $H_r^s(\RN)$.  In \cite{zjjjm}, the existence of radial ground state solutions was obtained when $V(x)\equiv V$ involving the critical growth by means of the constraint variational argument, where $V>0$ is a constant. When $V(x)=V(|x|)$, Secchi \cite{simonesecchi}  proved the existence of radially symmetric solutions for equation (\ref{p1}) in $H_r^s(\RN)$ by the fractional Pohoz\v{a}ev identity and the monotonicity trick in subcritical case. The conditions on $f$ and $V(x)$ are as follows
\begin{itemize}
\item [$(f_1')$] $f\in C(\R,\R)$ is of class $C^{1,\g}$ for some $\g>\max\{0,1-2s\}$, and odd,

\item [$(f_2')$]  $-\iy<\liminf_{t\to 0^+}{f(t)/t}\leq \limsup_{t\to 0^+}{f(t)/t}=-m<0$,

\item [$(f_3')$]  $-\iy<\limsup_{t\to +\iy}{f(t)/{t^{2_s^*-1}}}\leq 0$, where $2_s^*=\frac{2N}{N-2s}$,

\item [$(f_4')$] for some $\zeta>0$, there results $F(\zeta)=\int_0^\zeta f(t)dt>0$,

\item [$(H_1)$] $V\in C^1(\RN,\R),V(x)\geq 0$ for every $x\in \RN$ and this inequality is strict at some point,
\item [$(H_2)$] $\|\max\{\lan\na V(x),x\ran,0\}\|_{L^{N/2s}(\RN)}<2sS_s$,
\item [$(H_3)$] $\lim_{|x|\to +\iy}V(x)=0$,
\item [$(H_4)$] $V(x)$ is radially symmetric,
\end{itemize}
where $S_s$ is the best Sobolev constant for the critical embedding, that is
\begin{align*}
S_s=\inf_{u\in H^s(\RN),{u\neq 0}}\frac{\|(-\DD)^{s/2}u\|^2_{L^2}}{\|u\|^2_{L^{2_s^*}}},
\end{align*}
here $H^s(\RN)$ is the fractional Sobolev space with respect to the norm
\begin{align*}
\|u\|^2=\int_{\RN}|(-\DD)^{s/2}u|^2+|u|^2.
\end{align*}
Where $(f_2')-(f_4')$ are called Berestycki-Lions type conditions, and $(H_2)$ is used to get the boundedness of the PS sequence by use of the  monotonicity trick. The condition $f\in C^1$ ensure that the fractional Pohoz\v{a}ev identity can be used.

Now, the problem is how about the existence of ground state solutions when $V(x)$ is non-radially in the critical case. As we all know, for the critical case, the loss of the compactness for the embedding $H^s(\mathbb{R}^N)\hookrightarrow L^{2_s^*}(\mathbb{R}^N)$ is the main difficulty. What's more, PS condition, in general, fails. Since $V(x)$ is non-radial, the method introduced in \cite{zjjjm} can not be used here since the fractional space they used is $H_r^s(\RN)$. With critical growth, the authors \cite{xiaominghe} proved the existence of solutions for equation (\ref{p1}) under the monotonicity condition on $f(t)/t$  and  $0<\mu F(t)=\mu \int_0^tf(t)dt\leq tf(t),\mu\in(2,2_s^*)$ for all $t>0$.

Motivated by the seminal papers above, we use the indirect approach developed in \cite{L.Jeanjean} to investigate the existence of ground state solutions for nonlinear  fractional Schr\"{o}dinger equation (\ref{p1}) involving the critical nonlinearity, where the potential $V(x)$ depends on $x$ non-radially.
More precisely, on nonlinearity $f$, we assume
\begin{itemize}

\item [$(f_1)$] $f\in C^1(\R^+,\R)$ and $\lim_{t\to 0}f(t)/t=0$,

\item [$({f_2})$]  $\lim_{t\to \iy}{f(t)/t^{2_s^*-1}}=1$,

\item [$({f_3})$]  There exists $D>0$ and $p<2_s^*$ such that $f(t)\geq t^{2_s^*-1}+Dt^{p-1},t\geq 0$.

\end{itemize}
We assume $f(t)\equiv 0$ for $t\leq 0$  throughout the paper since we are concerned with the positive solutions.

On potential $V\in C^1(\RN,\R)$, we assume
\begin{itemize}
\item [$(V_1)$] There exists $V_0>0$ such that $\inf_{x\in \R ^N}V(x)\geq V_0$,
\item [$(V_2)$] $V(x)\leq V(\iy):=\lim_{|x|\to \iy}V(x)<\iy$ for all $x\in \RN$ and $V(x)\not \equiv  V(\iy)$,
\item [$(V_3)$] $\|\max\{\lan\na V(x),x\ran,0\}\|_{L^{N/2s}(\RN)}<2sS_s$.
\end{itemize}

In contrast to the conditions in \cite{xiaominghe}, our conditions are more weaker.
The main result is  the following.
\bt \lab{Theorem 1.1}
Assume $N>2s,s\in(0,1)$, if $\max\{2,2_s^*-2\}<p<2_s^*$, $(V_1)-(V_3)$ and $(f_1)-(f_3)$ hold, then problem $(\ref{p1})$ has a ground state solution.
\et

The proof of Theorem \ref{Theorem 1.1} is inspired by the ideas in \cite{L.JeanjeanandKazunagaTanaka} and \cite{J.Zhang}.

First, we show the existence of positive solutions of $(\ref{p1})$. For this purpose, we look for a special BPS sequence for the energy functional $I$ associated with $(\ref{p1})$ by use of the Struwe's monotonicity trick. Precisely, with the help of auxiliary energy functional $I_\la$ satisfying
\begin{align*}
I(u_{\la_j})=I_{\la_j}(u_{\la_j})+(\la_j-1)\int_{\RN} F(u_{\la_j}),\la_j\to 1,\ \ j\to \iy,
\end{align*}
we prove the existence of positive critical points denoted by $u_{\la_j}$ of $I_{\la_j}$. Thanks to the decomposition of BPS sequence, the properties of $\{u_{\la_j}\}$ and the energy estimation of $I_{\la_j}(u_{\la_j})$ are obtained. Consequently, we show that $\{u_{\la_j}\}$ is a BPS sequence for $I$ at some level value.

Secondly, for the proof of the existence of ground state solutions, we construct a minimizing sequence $\{u_n\}$ which is composed of the critical points of $I$. We show that $\{u_n\}$ is a BPS sequence for $I$ at $m$, here $m$ denote the least energy. Then, making use of the decomposition of BPS sequence and the relationship of $I$ and $I^\iy$, we prove that $m$ is attained at some $\ti u\neq 0$.

\br In the proof of our main results, the estimations of the Mountain Pass (MP for short) values, Pohoz\v{a}ev identity and the decomposition of BPS all play crucial roles.
\er
The paper is organized as follows. In section 2, we introduce a variational setting of our problem and present some preliminary results. In section 3, we are concerned with the decomposition of BPS and the existence of nontrivial critical points for the auxiliary energy functional. Section 4 is devoted to the completion of the proof of Theorem \ref{Theorem 1.1}.

{\bf Note}: in the following, the letters $C,\delta, \delta_0$ are indiscriminately used to denote various positive constants whose exact values are irrelevant.

\s {Preliminaries and functional setting}
\renewcommand{\theequation}{2.\arabic{equation}}
In order to establish the variational setting for (\ref{p1}), we give some useful facts of the fractional Sobolev space \cite{Nezza} and some preliminary lemmas.

The fractional Laplacian operator $(-\Delta)^s$ with $s\in(0,1)$ of a function $u:\mathbb{R}^N\rightarrow \mathbb{R}$ is defined by
\begin{align*}
\mathcal{F}((-\DD)^s u)(\xi)=|\xi|^{2s}\mathcal{F}(u)(\xi),\ \ \xi\in\mathbb{R}^N,
\end{align*}
where $\mathcal{F}$ is the Fourier transform. For $s\in (0,1)$, the fractional order Sobolev space $H^s(\mathbb{R}^N)$ is defined by
\begin{align*}
H^s(\mathbb{R}^N)=\{u\in L^2(\mathbb{R}^N):\int_{\mathbb{R}^N}|\xi|^{2s}|\hat{u}|^2d\xi<\infty\},
\end{align*}
endowed with the norm $\|u\|_{H^s(\mathbb{R}^N)}=(\int_{\mathbb{R}^N}(|\xi|^{2s}|\hat{u}|^2+|\hat{u}|^2)d\xi)^{\frac{1}2}$, where $\hat{u}\doteq \mathcal{F}(u)$.
By Plancherel's theorem, we have $\|u\|_{L^2(\mathbb{R}^N)}=\|\hat{u}\|_{L^2(\mathbb{R}^N)}$ and
\begin{align*}
\int_{\mathbb{R}^N}|(-\Delta)^{\frac{s}2}u(x)|^2dx
=\int_{\mathbb{R}^N}(|\xi|^s|\hat{u}|)^2d\xi
=\int_{\mathbb{R}^N}|\xi|^{2s}|\hat{u}|^2d\xi.
\end{align*}
It follows that $\|u\|_{H^s(\mathbb{R}^N)}=(\int_{\mathbb{R}^N}(|(-\Delta)^{\frac{s}2}u(x)|^2+|u|^2)dx)^{\frac{1}2}$, $u\in H^s(\mathbb{R}^N)$.
If $u$ is smooth enough, $(-\DD)^s u$ can be computed by the following singular integral
\begin{align*}
(-\DD)^s u(x)=c_{N,s}\mbox{P.V.}\int_{\RN}\frac{u(x)-u(y)}{|x-y|^{N+2s}}dy.
\end{align*}
Here $c_{N,s}$ is the normalization constant and P.V. is the principal value. So, one can get an alternative definition of the fractional Sobolev space $H^s(\mathbb{R}^N)$ as follows,
\begin{align*}
H^s(\mathbb{R}^N)=\{u\in L^2(\RN):\frac{|u(x)-u(y)|}{|x-y|^{\frac{N+2s}{2}}}\in L^2(\RN\times\RN)\}£¬
\end{align*}
with the norm
\begin{align*}
\|u\|_{H^s(\mathbb{R}^N)}=\left(\int_{\RN}|u|^2+\int_{{\R}^{2N}}\frac{|u(x)-u(y)|^2}{|x-y|^{N+2s}}\right)^{\frac{1}2}.
\end{align*}

The space
$D^s(\mathbb{R}^N)$ denotes the completion of $C_0^\infty(\RN)$ with respect to the Gagliardo norm
\begin{align*}
\|u\|_{D^s(\mathbb{R}^N)}=\left(\int_{\mathbb{R}^N}|\xi|^{2s}|\hat{u}|^2d\xi\right)^{\frac{1}2}=\left(\int_{\mathbb{R}^N}|(-\Delta)^{\frac{s}2}u|^2\right)^{\frac{1}2}.
\end{align*}

\noindent Since we investigate the existence of solutions of  problem  (\ref{p1}), we need the fractional Sobolev space  $H_{V}^s(\mathbb{R}^N)$ which is a Hilbert subspace of $H^s(\mathbb{R}^N)$ with the norm
\be \lab{norm}
\|u\|_{H_{V}^s(\mathbb{R}^N)}:=\left(\int_{\mathbb{R}^N}\left(|(-\Delta)^{\frac{s}2}u|^2+V(x)|u|^2\right)dx\right)^{\frac{1}2}<\infty.
\ee
It is easy to check that $H_{V}^s(\mathbb{R}^N)\equiv H^s(\mathbb{R}^N)$ if $(V1)-(V2)$ hold. In our paper, we shall work on $H^s(\mathbb{R}^N)$ with norm (\ref{norm})
and we denote $\|u\|_{H^s(\RN)}$ by $\|u\|$ for simplicity.

Associated to problem (\ref{p1}), the energy functional $I:H^s(\RN)\to \R$ is
\begin{align*}
I(u)=\frac{1}2\int_{\RN}|(-\DD)^{s/2}u|^2+V(x)|u|^2-\int_{\RN}F(u),u\in H^s(\RN),
\end{align*}
where $F(u)=\int_0^uf(t)dt$. The conditions $(f_1)-(f_3)$ imply that $I\in C^1(H^s(\RN),\R)$.

\bd
$u$ is said to be a solution of $(\ref{p1})$ if $u$ is a critical point of the energy functional $I$ and satisfies
\begin{align*}
\int_{\RN}(-\DD)^{\frac{s}{2}}u(-\DD)^{\frac{s}{2}}\vp+\int_{\RN}V(x)u\vp=\int_{\RN}f(u)\vp,u\in H^s(\RN),\forall \vp\in C_0^\infty(\RN).
\end{align*}
$u$ is said to be a ground state solution of $(\ref{p1})$ if $u$ is a solution with the least energy among all nontrivial solutions of $(\ref{p1})$.
\ed

In the paper, we use the embedding lemma and  Lions lemma  as follows.

\bl\lab{lemma1.1}(\cite{lionsembedding})
For any $s\in(0,1)$, $H^s(\RN)$ is continuously embedding into $L^r(\RN)$ for $r\in [2,2_s^*]$ and compactly embedding into $L_{loc}^r(\RN)$ for $r\in [2,2_s^*)$.
\el

\bl\lab{lionslemma}(\cite{lionslemma})
Suppose that $\{u_n\}$ is bounded in $H^s(\RN)$ and
\begin{align*}
\sup_{z\in \RN}\int_{B_1(z)}|u_n|^2\to 0.
\end{align*}
Then $\|u_n\|_{L^r}\to 0$ for $r\in (2,2_s^*)$ when $N\geq3$ and for $r\in (2,+\iy)$ when $N=1,2$. Here $B_1(z)=\{y\in\RN,|y-z|\leq1\}.$
\el

\s {Solutions for auxiliary problems}
\renewcommand{\theequation}{3.\arabic{equation}}
In this section, we consider the family of functionals $I_\la(u):H^s(\RN)\to \R$ defined by
\begin{align*}
I_\la(u)=\frac{1}2\int_{\RN}|(-\DD)^{s/2}u|^2+V(x)|u|^2-\la\int_{\RN}F(u).
\end{align*}
The corresponding auxiliary problems are
\be\lab{auxiliaryproblem}
(-\DD)^su+V(x)u=\la f(u).
\ee
The main aim of this section is to prove that for almost every $\la\in[\frac{1}2,1]$, $I_\la$ has a nontrivial critical point $u_\la$ such that $I_\la(u_\la)\leq c_\la$, where
\begin{align*}
c_\la=\inf_{\g\in\G}\max_{t\in[0,1]}I_\la(\g(t)),
\end{align*}
and
\begin{align*}
\G=\{\g\in C([0,1],H^s(\RN)),\g(0)=0 \ \ \mbox{and}\ \ I_\la(\g(1))<0\}.
\end{align*}

Before we proving the existence of solutions for the auxiliary problems (\ref{auxiliaryproblem}), we give some propositions and lemmas.
\bo \lab{pohozaev identity}
Let $u(x)$ be a critical point of $I_\la$ with $\la\in[\frac{1}2,1]$, then $u(x)$ satisfies
\be \lab{pohozaev identity}
\frac{N-2s}2\int_{\RN}|(-\DD)^{\frac{s}2}u|^2+\frac{N}2\int_{\RN}V(x)|u|^2+\frac{1}2\int_{\RN}\lan \na V(x),x\ran |u|^2-N\la\int_{\RN}F(u)=0.
\ee
\eo
As we all know, (\ref{pohozaev identity}) is named Pohoz\v{a}ev identity. The proof is similar as that in \cite{changxiaojun2} and we omit it here.

\bl \lab{splittinglemma} Assume $(f_1)-(f_2)$ hold. Let $\{u_n\}\subset H^s(\RN)$ be such that $u_n\to u$ weakly in $H^s(\RN)$. Then up to a subsequence,
\begin{align*}
\int_{\RN}(f(u_n)-f(u)-f(u_n-u))\phi=o_n(1)\|\phi\|.
\end{align*}
where $o_n(1)\to 0$ uniformly for $\phi\in C_0^\iy(\RN)$ as $n\to \iy$.
\el
\bp
The proof of this splitting lemma in the critical case on the fractional problems is  similar as that in \cite{zjjzwm}. So we omit the details.
\ep

Similar the proof of Brezis-Lieb Lemma in \cite{brezis-lieb}, we can give the following lemma.

\bl\lab{Brezis-Lieb lemma} For $s\in (0,1)$, assume $(f_1)-(f_2)$. Let $\{u_n\}\subset H^s(\RN)$ such that $u_n\to u$ weakly in $H^s(\RN)$ and a.e. in $\RN$ as $n\to \iy$, then
\begin{align*}
\int_{\RN}F(u_n)=\int_{\RN}F(u_n-u)+\int_{\RN}F(u)+o_n(1),
\end{align*}
where $o_n(1)\to 0$ as $n\to \iy$.
\el

In order to obtain the existence of critical points for $I_\la$, the following abstract result in \cite{L.Jeanjean} is needed, which shows that for almost every  $\la\in[\frac{1}2,1]$, $I_\la$ possesses a BPS  sequence at the level $c_\la$.

\bt \lab{Theorem 2.1} Let $X$ be a Banach space equipped with a norm $\|\cdot\|_X$ and let $J\subset\R^+$ be an interval. For a family $(I_\la)_{\la\in J}$ of $C^1$-functionals on $X$ of the form
\begin{align*}
I_\la(u)=A(u)-\la B(u),\forall \la\in J,
\end{align*}
where $B(u)\geq 0,\forall u\in X$ and such that either $A(u)\to +\iy$ or $B(u)\to +\iy$ as $\|u\|_X\to \iy$.
If there are two points $v_1,v_2$ in $X$ such that
\begin{align*}
c_\la=\inf_{\g\in\G}\max_{t\in[0,1]}I_\la(\g(t))>\max\{I_\la(v_1),I_\la(v_2)\},\forall \la\in J,
\end{align*}
where
\begin{align*}
\G=\{\g\in C([0,1],X),\g(0)=v_1,\g(1)=v_2\}.
\end{align*}
Then, for almost every $\la\in J$, there is a sequence $\{v_n\}\subset X$ such that
\begin{align*}
(i)\{v_n\} \ \ is\ \ bounded,\ \ (ii)I_\la(v_n)\to c_\la,\ \ (iii)I'_\la(v_n)\to 0\ \ in\ \ the\ \  dual\ \ X^{-1}\ \ of\ \ X.
\end{align*}
\et

In the following, we use Theorem \ref{Theorem 2.1} to seek nontrival critical points of $I_\la$ for almost every $\la\in J$. In what follows, let $X=H^s(\RN)$ and
\begin{align*}
A(u)=\frac{1}2\int_{\RN}|(-\DD)^{s/2}u|^2+V(x)|u|^2, B(u)=\int_{\RN}F(u).
\end{align*}
Obviously, $A(u)\to +\iy$ as $\|u\|\to \iy$ and $B(u)\geq 0$ for any $u\in H^s(\RN)$ by $(f_3)$. Now, we give the following lemma to ensure that $I_\la$ has MP geometry. Consequently, we get a BPS for $I_\la$ by Theorem \ref{Theorem 2.1}.

\bl \lab{Lemma 2.2}
Assume $(f_1)-(f_3),(V_1)-(V_2)$ hold. Then,

 (i)there exists a $v\in H^s(\RN)\setminus\{0\}$ with $I_\la(v)\leq 0$ for all $\la\in [\frac{1}2,1]$.

(ii)$c_\la=\inf_{\g\in\G}\max_{t\in[0,1]}I_\la(\g(t))>\max\{I_\la(0),I_\la(v)\}>0$ for all $\la\in [\frac{1}2,1]$. \\
         Here
\begin{align*}
\G=\{\g\in C([0,1],H^s(\RN)),\g(0)=0,\g(1)=v\}.
\end{align*}

(iii) there exists a BPS sequence $\{u_n\}$ at the MP level $c_\la$ for $I_\la$. Here $u_n\geq 0$.
\el

\bp If $(f_1)-(f_2)$ hold, then for any $\e>0$, there exists $C(\e)>0$ such that
\begin{align*}
\int_{\RN}F(u)\leq \e\int_{\RN}|u|^2+C(\e)\int_{\RN}|u|^{2_s^*},\forall u\in H^s(\RN).
\end{align*}
Thus
\begin{align*}
I_\la(u)&=\frac{1}2\|u\|^2-\la\int_{\RN}F(u)\\
&\geq \frac{1}2\|u\|^2-\e\|u\|_{L^2}^2-C(\e)\|u\|_{L^{2_s^*}}^{2_s^*}
\end{align*}
From Lemma \ref{lemma1.1}, there exist constant $\rho>0$  and $\delta>0$ independent of $\la$ such that for $\|u\|=\rho$, $I_\la(u)\geq \delta.$
On the other hand, $(f_3)$ implies that
\begin{align*}
I_\la(u)\leq\frac{1}2\|u\|^2-\frac{1}2\|u\|_{L^{2_s^*}}^{2_s^*}-\frac{D}{2p}\|u\|_{L^p}^p.
\end{align*}
Set $v_0\in H^s(\RN)$ such that $v_0\geq 0,v_0\neq 0$. Since $I_\la(tv_0)\to -\iy$ as $t\to+\iy$, then there exists $t_0$ such that $I_\la(t_0v_0)<0$ as $\|t_0v_0\|>\rho.$ Set $v=t_0v_0$, then $(i)$ and $(ii)$ hold. So, the conditions of Theorem \ref{Theorem 2.1} are satisfied. Therefore, for almost every $\la\in [\frac{1}2,1]$, there exists a BPS sequence $\{u_n\}$ for $I_\la$ at the MP value $c_\la$.
\noindent Now, we show $u_n\geq 0$. Let $u_n=u_n^++u_n^-$. Using $u_n^-$ as a test function, since $f(t)\equiv 0$ for all $t\leq 0$, we have
\begin{align*}
(I'_\la(u_n),u_n^-)&=\int_{\RN}(-\DD)^{\frac{s}2}u_n(-\DD)^{\frac{s}2}u_n^-+\int_{\RN}V(x)(u_nu_n^-)-\la\int_{\RN}f(u_n)u_n^-\\
&=\int_{\RN}(-\DD)^{\frac{s}2}u_n(-\DD)^{\frac{s}2}u_n^-+\int_{\RN}V(x)|u_n^-|^2.
\end{align*}
Since for every $x,y\in\RN$, we always have $\left(u_n^+(x)-u_n^+(y)\right)\left(u_n^-(x)-u_n^-(y)\right)\geq 0$, then
\begin{align*}
\left(u_n(x)-u_n(y)\right)\left(u_n^-(x)-u_n^-(y)\right)&=\left(u_n^+(x)-u_n^+(y)\right)\left(u_n^-(x)-u_n^-(y)\right)+\left(u_n^-(x)-u_n^-(y)\right)^2\\
&\geq \left(u_n^-(x)-u_n^-(y)\right)^2.
\end{align*}
Thus
\begin{align*}
\int_{\RN}(-\DD)^{\frac{s}2}u_n(-\DD)^{\frac{s}2}u_n^-&=\int_{\RN}\int_{\RN}\frac{\left(u_n(x)-u_n(y)\right)\left(u_n^-(x)-u_n^-(y)\right)}{|x-y|^{N+2s}}\\
&\geq \int_{\RN}\int_{\RN}\frac{\left(u_n^-(x)-u_n^-(y)\right)^2}{|x-y|^{N+2s}}\\
&=\int_{\RN}|(-\DD)^{\frac{s}2}u_n^-|^2.
\end{align*}
Therefore, from $(I'_\la(u_n),u_n^-)\to 0$, we have $\|u_n^-\|\to 0$. The proof is finished.
\ep

From the argument above, we obtain a BPS for $I_\la$ at the level $c_\la$. In order to get the convergence of the BPS sequence, we give some lemmas and propositions.

\bl \lab{Lemma 2.3}
Assume $(V_1)-(V_2)$ and $(f_1)-(f_3)$ hold, if $\max\{2,2_s^*-2\}<p<2_s^*$, then
\begin{align*}
c_\la<\frac{s}{N\la^{\frac{N-2s}{2s}}}S_s^{\frac{N}{2s}}.
\end{align*}
\el

\bp
Let $\vp\in C_0^\iy(\RN)$ is a cut-off function with support $B_2$ such that $\vp\equiv1$ on $B_1$ and $0\leq \vp\leq1$ on $B_2$, where $B_r$ denotes the ball in $\RN$ of center at origin and radius $r$. For $\e>0$, we define $\psi_\e(x)=\vp(x)U_\e(x)$, where
\begin{align*}
U_\e(x)=\kappa\e^{-\frac{N-2s}2}\left(\mu^2+\left|\frac{x}{\e S_s^{\frac{1}{2s}}}\right|^2\right)^{-\frac{N-2s}2}.
\end{align*}
By \cite{A.Cotsiolis}, $S_s$ can be achieved by $U_\e(x)$. Let $v_\e=\frac{\psi_\e}{\|\psi_\e\|_{L^{2_s^*}}}$, then $\|(-\DD)^{\frac{s}2}v_\e\|_{L^2}^2\leq S_s+O(\e^{N-2s})$. From \cite{xiaominghe}, we have the estimations,

\begin{align*}
\|v_\e\|_{L^2}^2=\left\{
\begin{array}{ll}
O(\e^{2s}),\ \ \ \ \ \ \ \ N>4s,\\
O(\e^{2s}\ln\frac{1}\e), \ \ \ N=4s,\\
O(\e^{N-2s}),\ \ \ \ \ N<4s,
\end{array}
\right.
\end{align*}
and

\begin{align*}
\|v_\e\|_{L^p}^p=\left\{
\begin{array}{ll}
O(\e^{\frac{2N-(N-2s)p}2}),\ \ \ \ p>\frac{N}{N-2s},\\
O(\e^{\frac{(N-2s)p}2}), \ \ \ \ \ \ \ \ \  p<\frac{N}{N-2s}.\\
\end{array}
\right.
\end{align*}

\noindent For any $t>0$, by $(f_3)$
\begin{align*}
I_\la(tv_\e)&=\frac{t^2}2\int_{\RN}|(-\DD)^{\frac{s}2}v_\e|^2+V(x)|v_\e|^2-\la\int_{\RN} F(tv_\e)\\
&=\frac{t^2}2\|v_\e\|^2-\la\int_{\RN} F(tv_\e)\\
&\leq \frac{t^2}2\|v_\e\|^2-\frac{\la}{2_s^*}t^{2_s^*}-\frac{Dt^p}{2p}\|v_\e\|_{L^p}^p.
\end{align*}
Obviously, $I_\la(tv_\e)\to -\iy$ as $t\to +\iy$ and $I_\la(tv_\e)>0$ for $t>0$
small.
\noindent Let $g(t)=\frac{t^2}2\|v_\e\|^2-\frac{\la}{2_s^*}t^{2_s^*}$. Then $t_\e=\left(\frac{\|v_\e\|^2}\la\right)^{\frac{1}{2_s^*-2}}$ is the maximum point of $g(t)$.

For $\e<1$, by the definition of $v_\e$, there exists $t_1>0$ small enough such that
\begin{align*}
\max_{t\in (0,t_1)}I_\la(tv_\e)\leq \frac{t^2}2\|v_\e\|^2<\frac{s}{N\la^{\frac{N-2s}{2s}}}S_s^{\frac{N}{2s}}.
\end{align*}

Since $I_\la(tv_\e)\to -\iy$ as $t\to +\iy$, it is easy to obtain that there exists $t_2>0$ such that
\begin{align*}
\max_{t\in (t_2,+\iy)}I_\la(tv_\e)<\frac{s}{N\la^{\frac{N-2s}{2s}}}S_s^{\frac{N}{2s}}.
\end{align*}

If $t\in [t_1,t_2]$,
\begin{align*}
\max_{t\in [t_1,t_2]}I_\la(tv_\e)&\leq \max_{t\in [t_1,t_2]}\{g(t)-\frac{Dt_1^p}{2p}\|v_\e\|_{L^p}^p\}  \\
&\leq g(t_\e)-\frac{Dt_1^p}{2p}\|v_\e\|_{L^p}^p\\
\end{align*}
For $g(t_\e)$, we have

\begin{align*}
g(t_\e)&=\frac{s}{N\la^{\frac{N-2s}{2s}}}(\|v_\e\|^2)^{\frac{N}{2s}}\\
&=\frac{s}{N\la^{\frac{N-2s}{2s}}}\left(\|(-\DD)^{\frac{s}2}v_\e\|_{L^2}^2+\int_{\RN}V(x)|v_\e|^2\right)^{\frac{N}{2s}}\\
&\leq \frac{s}{N\la^{\frac{N-2s}{2s}}}\left(S_s+O(\e^{N-2s})+C\|v_\e\|_{L^2}^2\right)^{\frac{N}{2s}}.
\end{align*}

\noindent By $(a+b)^q\leq a^q+q(a+b)^{q-1}b$, where $a>0,b>0,q>1$, we have
\begin{align*}
g(t_\e)&\leq \frac{s}{N\la^{\frac{N-2s}{2s}}}\left(S_s^{\frac{N}{2s}}+{\frac{N}{2s}}\left(S_s+O(\e^{N-2s})+C\|v_\e\|_{L^2}^2\right)^{\frac{N-2s}{2s}}\left(O(\e^{N-2s})+C\|v_\e\|_{L^2}^2\right)\right)\\
&\leq \frac{s}{N\la^{\frac{N-2s}{2s}}}S_s^{\frac{N}{2s}}+O(\e^{N-2s})+C\|v_\e\|_{L^2}^2.
\end{align*}

\noindent Thus
\begin{align*}
\max_{t\in [t_1,t_2]}I_\la(tv_\e)\leq \frac{s}{N\la^{\frac{N-2s}{2s}}}S_s^{\frac{N}{2s}}+O(\e^{N-2s})+C\|v_\e\|_{L^2}^2-\frac{Dt_1^p}{2p}\|v_\e\|_{L^p}^p.
\end{align*}

\noindent Next, we estimate $\max_{t\in [t_1,t_2]}I_\la(tv_\e)$ in three cases.

\noindent Case 1: If $N>4s$, then $\frac{N}{N-2s}<2$, together with $p>\max\{2,2_s^*-2\}$, we have $p>\frac{N}{N-2s}$. So
\begin{align*}
\max_{t\in [t_1,t_2]}I_\la(tv_\e)\leq \frac{s}{N\la^{\frac{N-2s}{2s}}}S_s^{\frac{N}{2s}}+O(\e^{N-2s})+O(\e^{2s})-O(\e^{\frac{2N-(N-2s)p}2}).
\end{align*}
From $p>2,N>4s$, then $\frac{2N-(N-2s)p}2<2s<N-2s$. Thus, for $\e>0$ small enough, we obtain

\begin{align*}
\max_{t\in [t_1,t_2]}I_\la(tv_\e)< \frac{s}{N\la^{\frac{N-2s}{2s}}}S_s^{\frac{N}{2s}}.
\end{align*}

\noindent Case 2: If $N=4s$, then $2<p<4$. For $\e>0$ small enough, we obtain
\begin{align*}
\max_{t\in [t_1,t_2]}I_\la(tv_\e)&\leq \frac{s}{N\la^{\frac{N-2s}{2s}}}S_s^{\frac{N}{2s}}+O(\e^{N-2s})+O(\e^{2s}\ln\frac{1}\e)-O(\e^{4s-sp})\\
&\leq\frac{s}{N\la^{\frac{N-2s}{2s}}}S_s^{\frac{N}{2s}}+O\left(\e^{2s}(1+\ln\frac{1}\e)\right)-O(\e^{4s-sp})\\
&<\frac{s}{N\la^{\frac{N-2s}{2s}}}S_s^{\frac{N}{2s}},
\end{align*}
since
\begin{align*}
\lim_{\e\to 0^+}\frac{\e^{4s-sp}}{\e^{2s}(1+\ln\frac{1}\e)}\to +\iy.
\end{align*}

\noindent Case 3: If $2s<N<4s$, then $\frac{N}{N-2s}>2$, together with $p>\max\{2,2_s^*-2\}$, we have $p>\frac{N}{N-2s}$. So
\begin{align*}
\max_{t\in [t_1,t_2]}I_\la(tv_\e)\leq \frac{s}{N\la^{\frac{N-2s}{2s}}}S_s^{\frac{N}{2s}}+O(\e^{N-2s})-O(\e^{\frac{2N-(N-2s)p}2}).
\end{align*}
From $p>\frac{4s}{N-2s}$, then $\frac{2N-(N-2s)p}2<N-2s$. For $\e>0$ small enough, we obtain
\begin{align*}
\max_{t\in [t_1,t_2]}I_\la(tv_\e)< \frac{s}{N\la^{\frac{N-2s}{2s}}}S_s^{\frac{N}{2s}}.
\end{align*}
The proof is completed.
\ep

In (\ref{p1}), if $V(x)\equiv V(\iy)$, for $\la \in [\frac{1}2,1]$, the family of functionals $I_\la^\iy:H^s(\RN)\mapsto \R$, defined as
\begin{align*}
I_\la^\iy(u)=\frac{1}2\int_{\RN}|(-\DD)^{\frac{s}2}u|^2+V(\iy)|u|^2-\la\int_{\RN}F(u),
\end{align*}
plays an important role in our paper.
Similar  as that in \cite{JeanTa,zjjzwm2}, we can derive the following result.
\bl \lab{Lemma 2.5} For $\la \in [\frac{1}2,1]$, if $w_\la\in H^s(\RN)$ is a nontrivial critical point of $I_\la^\iy$, then there exists $\g_\la\in C([0,1],H^s(\RN))$ such that $\g_\la(0)=0,I_\la^\iy(\g_\la(1))<0,w_\la\in \g_\la[0,1]$ and $\max_{t\in[0,1]}I_\la^\iy(\g_\la(t))=I_\la^\iy(w_\la)$.
\el

\bl \lab{Lemma 2.6}\cite{zjjjm} If $f$ satisfies $(f_1)-(f_3),\max\{2,2_s^*-2\}<p<2_s^*$, then for almost every $\la \in [\frac{1}2,1]$, $I_\la^\iy$ has a positive ground state solution.
\el

\bl\lab{energyest}
If $V(x)\equiv V(\iy)>0$ and $(f_1)-(f_2)$ hold, there exists a constant $\delta> 0$ independent of $\la$ such that any nontrivial critical point $u$ of $I_\la^\iy$ satisfies $I_\la^\iy(u)\geq\delta$.
\bp
Letting $u$ is a nontrivial critical point of $I_\la^\iy$, from Pohoz\u{a}ev identity (\ref{pohozaev identity}), we have
\begin{align*}
I_\la^\iy(u)=\frac{s}N\int_{\RN}|(-\DD)^{\frac{s}2}u|^2.
\end{align*}
If $(f_1)-(f_2)$ hold, then for any $\e>0$, there exists $C(\e)>0$ such that
\begin{align*}
\int_{\RN}|(-\DD)^{\frac{s}2}u|^2+V|u|^2\leq \e\int_{\RN}|u|^2+C(\e)\int_{\RN}|u|^{2_s^*}.
\end{align*}
Thus, $\int_{\RN}|(-\DD)^{\frac{s}2}u|^2\leq C\int_{\RN}|u|^{2_s^*}$. On the other hand, by the Sobolev embedding theorem, we have $\int_{\RN}|u|^{2_s^*}\leq\ti C (\int_{\RN}|(-\DD)^{\frac{s}2}u|^2)^{\frac{2_s*}2}$. Since $u\neq 0$, there exists a constant $\delta_0 >0$ such that $\int_{\RN}|(-\DD)^{\frac{s}2}u|^2\geq \delta_0$ and so $I_\la^\iy(u)\geq \delta:=s\delta_0/N$. The proof is finished.
\ep
\el

Now, we  give the decomposition of a BPS sequence.

\bo \lab{proposition1} Assume $(V_1)-(V_3)$ and $(f_1)-(f_3)$ hold, if $\max\{2,2_s^*-2\}<p<2_s^*$,  for almost every $\la\in [\frac{1}2,1]$, $\{u_n\}$ given in Lemma \ref{Lemma 2.2} is the BPS sequence at the MP value $c_\la$. Moreover,
$c_\la<\frac{s}{N\la^{\frac{N-2s}{2s}}}S_s^{\frac{N}{2s}}.$
Then there exist a subsequence, still denoted by $\{u_n\}$, an integer $k\in \NN\cup\{0\}$ and $v_\la^j\in H^s(\RN)$ for $1\leq j\leq k$, such that
\begin{itemize}
\item [$(i)$] $u_n\to u_\la$ weakly in $H^s(\RN)$ and $I_\la'(u_\la)=0$,
\item [$(ii)$] $v_\la^j\neq0,v_\la^j\geq 0$ and $I_\la^{\iy'}(v_\la^j)=0$ for $1\leq j\leq k$,
\item [$(iii)$] $c_\la=I_\la(u_\la)+\sum_{j=1}^k I_\la^\iy(v_\la^j)$,
\item [$(iv)$] $\|u_n-u_0-\sum_{j=1}^k v_\la^j(\cdot-y_n^j)\|\to 0$.
\end{itemize}
where  $|y_n^j|\to\iy$ and $|y_n^i-y_n^j|\to \iy$ as $n\to \iy$ for any $i\neq j$.
\eo

\bp For $\la\in [\frac{1}2,1]$, let $\{u_n\}\subset H^s(\RN),u_n\geq0$ be given in Lemma \ref{Lemma 2.2}. Since $\{u_n\}$ is bounded, there exist a subsequence denoted by $\{u_n\}$ and $u_\la\in H^s(\RN)$ satisfying $u_n\to u_\la$ weakly  in $H^s(\RN)$ and $u_n\to u_\la$ a.e. in $\RN$. It is not hard to verify that $I_\la'(u_\la)=0$.

{\bf Step 1.} Let $v_n^1=u_n-u_\la$, if $v_n^1\to 0$ strongly in $H^s(\RN)$, the Proposition holds with $k=0$.

{\bf Step 2.} We  claim that if $v_n^1\nrightarrow 0$ strongly, $\lim_{n\to \iy}\sup_{z\in \RN}\int_{B_1(z)}|v_n^1|^2>0.$

Since $I_\la(u_n)\to c_\la$, by Lemma \ref{Brezis-Lieb lemma}, we have
\be\lab{guji1}
c_\la-I_\la(u_\la)=I_\la(v_n^1)+o(1).
\ee

Since $v_n^1\rightharpoonup 0$, by $(V_2)$ and Lemma \ref{lemma1.1}, we have
\begin{align*}
I_\la^\iy(v_n^1)-I_\la(v_n^1)&=\int_{\RN}(V(\iy)-V(x))|v_n^1|^2\\
&=\int_{B_{R}(0)}(V(\iy)-V(x))|v_n^1|^2+\int_{\RN\backslash B_{R}(0)}(V(\iy)-V(x))|v_n^1|^2\\
&\to 0.
\end{align*}
Consequently,
\be\lab{guji2}
c_\la-I_\la(u_\la)=I_\la^\iy(v_n^1)+o(1).
\ee

Suppose $\lim_{n\to \iy}\sup_{z\in \RN}\int_{B_1(z)}|v_n^1|^2=0.$ By  Lemma \ref{lionslemma}, we have
\be
\label{lions1}
v_n^1\to 0 \ \ \mbox{in} \ \ L^t(\RN), \ \ \forall t\in(2,2_s^*).
\ee
Let $f(t)=h(t)+(t^+)^{2_s^*-1}$, from $(f_1)-(f_2)$, for any $\e>0$, there exists $C(\e)>0$ such that
\begin{align*}
\left|\int_{\RN}H(v_n^1)\right|\leq \e\left(\int_{\RN}|v_n^1|^2+|v_n^1|^{2_s^*}\right)+C(\e)\int_{\RN}|v_n^1|^r,
\end{align*}
where $r<2_s^*$. Since $v_n^1\in H^s(\RN)$, together with (\ref{lions1}), we obtain
\begin{align*}
\left|\int_{\RN}H(v_n^1)\right|\leq \e C+o(1),
\end{align*}
which implies $\int_{\RN}H(v_n^1)=o(1)$ since $\e$ is small enough.
What's more, by Brezis-Lieb lemma, we have
\begin{align*}
\int_{\RN}|v_n^1|^{2_s^*}=\int_{\RN}|u_n|^{2_s^*}-\int_{\RN}|u_\la|^{2_s^*}+o(1).
\end{align*}
Thus, (\ref{guji1}) reduces to
\be\lab{guji3}
c_\la-I_\la(u_\la)=\frac{1}2\|v_n^1\|^2-\frac{\la}{2_s^*}\int_{\RN}|v_n^1|^{2_s^*}+o(1).
\ee
Noting that  $I'_\la(u_n)v_n^1\to 0$ and $I'_\la(u_\la)v_n^1=0$, by direct calculation, we get
\begin{align*}
\|v_n^1\|^2-\la\int_{\RN}\left(f(u_n)-f(u_\la)\right)v_n^1=I'_\la(u_n)v_n^1-I'_\la(u_\la)v_n^1\to 0.
\end{align*}
By Lemma \ref{splittinglemma},
\begin{align*}
\int_{\RN}(f(u_n)-f(u))v_n^1=\int_{\RN}f(v_n^1)v_n^1+o(1)\|v_n^1\|=\int_{\RN}h(v_n^1)v_n^1+\int_{\RN}|v_n^1|^{2_s^*}+o(1)\|v_n^1\|.
\end{align*}
By (\ref{lions1}) and similar argument as above, we have
$\int_{\RN}\left(f(u_n)-f(u_\la)\right)v_n^1=\int_{\RN}|v_n^1|^{2_s^*}+o(1)$. Therefore,
\be\lab{guji4}
\|v_n^1\|^2-\la\int_{\RN}|v_n^1|^{2_s^*}=o(1).
\ee
Combining (\ref{guji3}) with (\ref{guji4}), we get $c_\la-I_\la(u_\la)=\frac{s}N\|v_n^1\|^2+o(1)$.

Noting that $I_\la'(u_\la)=0$, from Pohoz\v{a}ev identity (\ref{pohozaev identity}) and Sobolev embedding theorem , we obtain
\begin{align*}
I_\la(u_\la)&=\frac{s}N\int_{\RN}|(-\DD)^{\frac{s}2}u_\la|^2-\frac{1}{2N}\int_{\RN}\lan\na V(x),x\ran u_\la^2\\
&\geq \frac{s}N\int_{\RN}|(-\DD)^{\frac{s}2}u_\la|^2-\frac{1}{2NS_s}\|\max\{\lan\na V(x),x\ran,0\}\|_{L^{\frac{N}{2s}}}\int_{\RN}|(-\DD)^{\frac{s}2}u_\la|^2.
\end{align*}
$(V_3)$ implies $I_\la(u_\la)\geq 0$.
Thus $c_\la-I_\la(u_\la)<\frac{s}{N\la^{\frac{N-2s}{2s}}}S_s^{\frac{N}{2s}}$. On the other hand, since $v_n^1\nrightarrow 0$ strongly,  there exists a constant $l>0$ such that $\|v_n^1\|^2\to l$. Set $\|(-\DD)^{\frac{s}2}v_n^1\|_{L^2}^2=\ti l<l$, then
\begin{align*}
S_s=\inf_{u\in H^s(\RN),{u\neq 0}}\frac{\|(-\DD)^{s/2}u\|^2_{L^2}}{\|u\|^2_{L^{2_s^*}}}\leq \frac{\ti l}{(\frac{l}\la)^{\frac{2}{2_s^*}}}\leq l^{\frac{2s}{N}}\la^{\frac{N-2s}N}.
\end{align*}
So we have $l\geq \frac{S_s^{\frac{N}{2s}}}{\la^{\frac{N-2s}{2s}}}$. Consequently, $c_\la-I_\la(u_\la)\geq\frac{s}{N\la^{\frac{N-2s}{2s}}}S_s^{\frac{N}{2s}}$, which is a contradiction. The claim is true.

{\bf Step 3.} From the argument in step 2, if $v_n\rightharpoonup 0$, then  $\lim_{n\to \iy}\sup_{z\in \RN}\int_{B_1(z)}|v_n^1|^2>0.$
Thus, after extracting a subsequence if necessary, there exist $\{z_n^1\}\subset \RN$ and  $v_\la^1\in H^s(\RN)$ such that $|z_n^1|\to\iy$ and
\begin{align*}
(i)\lim_{n\to \iy}\int_{B_1(z_n^1)}|v_n^1|^2>0, \ \ (ii)v_n^1(\cdot+z_n^1)\rightharpoonup v_\la^1\neq 0,\ \ (iii)I_\la^{\iy'} (v_\la^1)=0.
\end{align*}
Clearly (i),(ii) are standard and the point is to show (iii). Set $u_n^1=v_n^1(\cdot+z_n^1)$. To prove $I_\la^{\iy'} (v_\la^1)=0$, it suffices to prove  $I_\la^{\iy'} (u_n^1)\to 0$. For any $\vp\in C_0^\iy(\RN)$, from $I_\la'(v_n^1)\to 0$, we have
\begin{align*}
I_\la'(v_n^1)\vp(\cdot -z_n^1)&=\int_{\RN}(-\DD)^{\frac{s}2}v_n^1(x+z_n^1)(-\DD)^{\frac{s}2}\vp(x)dx+\int_{\RN}V(x+z_n^1)v_n^1(x+z_n^1)\vp(x)dx\\
&-\int_{\RN}f(v_n^1(x+z_n^1))\vp(x)dx\\
&=\int_{\RN}(-\DD)^{\frac{s}2}u_n^1(x)(-\DD)^{\frac{s}2}\vp(x)dx+\int_{\RN}V(x+z_n^1)u_n^1(x)\vp(x)dx\\
&-\int_{\RN}f(u_n^1(x))\vp(x)dx\\
&\to 0.
\end{align*}
Since $|z_n^1|\to\iy$ and $\vp\in C_0^\iy(\RN)$, by $(V_2)$, we get
\begin{align*}
\int_{\RN}V(x+z_n^1)u_n^1(x)\vp(x)dx\to \int_{\RN}V(\iy)u_n^1(x)\vp(x)dx.
\end{align*}
Thus,
\begin{align*}
I_\la^{\iy'}(u_n^1)=\int_{\RN}(-\DD)^{\frac{s}2}u_n^1(x)(-\DD)^{\frac{s}2}\vp(x)dx+\int_{\RN}V(\iy)u_n^1(x)\vp(x)dx-\int_{\RN}f(u_n^1(x))\vp(x)dx\to 0.
\end{align*}
Then we get $I_\la^{\iy'} (v_\la^1)=0$ since $u_n^1\rightharpoonup v_\la^1$. On the other hand, from (\ref{guji2}), it is easy to see that $c_\la-I_\la(u_\la)=I_\la^\iy(u_n^1)+o(1).$

So, we get a bounded sequence $\{u_n^1\}$ with $u_n^1\rhu v_\la^1\neq 0$ satisfying
\begin{align*}
I_\la^\iy(u_n^1)\to c_\la-I_\la(u_\la),\ \ \ I_\la^{\iy'}(u_n^1)\to 0,\ \ \ I_\la^{\iy'} (v_\la^1)=0.
\end{align*}
Let $v_n^2=u_n^1-v_\la^1$. Then $u_n=u_\la+v_\la^1(\cdot-z_n^1)+v_n^2(\cdot-z_n^1)$. If $v_n^2\to 0$ strongly in $H^s(\RN)$, we have
\begin{align*}
\left\{
\begin{array}{ll}
c_\la-I_\la(u_\la)=I_\la^\iy(v_\la^1),\\
\|u_n-u_\la-v_\la^1(\cdot-z_n^1)\|\to 0.
\end{array}
\right.
\end{align*}
If $v_n^2\nrightarrow 0$ strongly, similarly as (\ref{guji1}) and (\ref{guji2}), we have
\begin{align*}
c_\la-I_\la(u_\la)-I_\la^\iy(v_\la^1)=I_\la^\iy(v_n^2)+o(1),I_\la^{\iy'}(v_n^2)\to 0.
\end{align*}
By the same argument as step 2, we obtain $\lim_{n\to \iy}\sup_{z\in \RN}\int_{B_1(z)}|v_n^2|^2>0.$ Then, there exist $\{z_n^2\}\subset \RN$ and  $v_\la^2\neq 0$ such that $|z_n^2|\to\iy$ and
\begin{align*}
(i)\lim_{n\to \iy}\int_{B_1(z_n^2)}|v_n^1|^2>0, \ \ (ii)v_n^2(\cdot+z_n^2)\rightharpoonup v_\la^2,\ \ (iii)I_\la^{\iy'} (v_\la^2)=0.
\end{align*}
Set $u_n^2=v_n^2(\cdot+z_n^2)$. Then, $\{u_n^2\}$ is a bounded sequence satisfying $u_n^2\rhu v_\la^2$ and
\begin{align*}
I_\la^\iy(u_n^2)\to c_\la-I_\la(u_\la)-I_\la^\iy(v_\la^1),\ \ \ I_\la^{\iy'}(u_n^2)\to 0.
\end{align*}
Let $v_n^3=u_n^2-v_\la^2$. Then $u_n=u_\la+v_\la^1(\cdot-z_n^1)+v_\la^2(\cdot-z_n^1-z_n^2)+v_n^3(\cdot-z_n^1-z_n^2)$. If $v_n^3\to 0$ strongly in $H^s(\RN)$, we have
\begin{align*}
\left\{
\begin{array}{ll}
c_\la=I_\la(u_\la)+I_\la^\iy(v_\la^1)+I_\la^\iy(v_\la^2),\\
\|u_n-u_\la-v_\la^1(\cdot-z_n^1)-v_\la^2(\cdot-z_n^1-z_n^2)\|\to 0.
\end{array}
\right.
\end{align*}
Otherwise, we repeat the procedure above. From Lemma \ref{energyest}, we can terminate our arguments by repeating the above proof by finite $k$ steps. That is, let $y_n^j=\sum_{i=1}^{j}z_n^i$, then
\begin{align*}
\left\{
\begin{array}{ll}
c_\la=I_\la(u_\la)+\sum_{j=1}^{k}I_\la^\iy(v_\la^j),\\
\|u_n-u_\la-\sum_{j=1}^{k}v_\la^j(\cdot-y_n^j)\|\to 0.
\end{array}
\right.
\end{align*}
{\bf Step 4.} Now, we show that after extracting a subsequence of $\{y_n^j\}$ and redefining $\{v_\la^j\}$ if necessary, $(iii),(iv)$ hold for $|y_n^j|\to\iy$ and $|y_n^i-y_n^j|\to \iy$ as $n\to \iy$ for any $i\neq j$.
Let $A=\{1,2,\cdot\cdot\cdot,k\}$.  From $u_n-u_\la-\sum_{j=1}^{k}v_\la^j(\cdot-y_n^j)\to 0$ and $u_n\to u_\la$ a.e. in $\RN$, we get that $\sum_{j=1}^{k}v_\la^j(\cdot-y_n^j)\to 0$ a.e. in $\RN$. Since $v_\la^j\geq0$ for any $j$, then $|y_n^j|\to \iy$. For $y_n^i$, assume $A_i=\{y_n^j:|y_n^i-y_n^j|\ \mbox{is bounded for }n \}$, then up to a sequence, there exists some $\ti v_\la^i\in H^s(\RN)$ such that $\sum_{j\in A_i}v_\la^j(\cdot+y_n^i-y_n^j)\to \ti v_\la^i$ strongly in $H^s(\RN)$. Then $\|u_n-u_\la-\ti v_\la^i(\cdot-y_n^i)-\sum_{j\in (A\setminus A_i)}v_\la^j(\cdot-y_n^j)\|\to 0$. Since $v_\la^j(j\in A)$ is the critical point of $I_\la^{\iy'}$, we have $I_\la^{\iy'}(\ti v_\la^i)=0$. Then we redefine $v_\la^i:=\ti v_\la^i$, and then $\|u_n-u_\la-\sum_{j\in (A\setminus A_i)\cup \{i\}}v_\la^j(\cdot-y_n^j)\|\to 0$ holds as $n\to \iy$. By repeating the argument above at most $(k-1)$ times and redefining $\{v_\la^j\}$ if necessary, there exists $\Lambda\subset A$ such that
\begin{align*}
\left\{
\begin{array}{ll}
|y_n^j|\to\iy,|y_n^i-y_n^j|\to \iy,\forall i\neq j, n\to \iy, \\
\|u_n-u_\la-\sum_{j\in \Lambda}v_\la^j(\cdot-y_n^j)\|\to 0.
\end{array}
\right.
\end{align*}
The proof is finished.
\ep

If $V(x)\equiv V>0$, we can get the similar decomposition of the BPS sequence for the autonomous problem (\ref{p1}). Denote the energy functional of autonomous problem (\ref{p1}) and auxiliary energy functional by $J$ and $J_\la$ $(\la\in [\frac{1}2,1])$ respectively. Let $c_\la$ be the MP value for $J_\la$, then we have the following result.

\bc\lab{corollary} Assume $V(x)\equiv V>0$ and $(f_1)-(f_3)$ hold. For $\la\in[\frac{1}2,1]$, if $\{u_n\}\subset H^s(\RN)$ is a sequence such that $u_n\geq0,\|u_n\|<\iy, J_\la(u_n)\to c_\la$ and $J_\la'(u_n)\to 0$, what's more $c_\la<\frac{s}{N\la^{\frac{N-2s}{2s}}}S_s^{\frac{N}{2s}}$. Then there exist a subsequence of $\{u_n\}$, an integer $l\in \NN\cup\{0\}$ and $w_\la^j\in H^s(\RN)$ for $1\leq j\leq l$ such that
\begin{itemize}
\item [$(i)$] $u_n\to u_\la$ weakly in $H^s(\RN)$ with $J_\la'(u_\la)=0$,
\item [$(ii)$] $w_\la^j\neq0,w_\la^j\geq 0$ and $J_\la'(w_\la^j)=0$ for $1\leq j\leq l$,
\item [$(iii)$] $c_\la=J_\la(u_\la)+\sum_{j=1}^l J_\la(w_\la^j)$,
\item [$(iv)$] $\|u_n-u_0-\sum_{j=1}^l w_\la^j(\cdot-y_n^j)\|\to 0$,
\end{itemize}
where  $|y_n^j|\to\iy$ and $|y_n^i-y_n^j|\to \iy$ as $n\to \iy$ for any $i\neq j$.
\ec

The proof is similar as Proposition \ref{proposition1}, we omit it here.

Now, we complete the proof of the existence of solutions of the auxiliary problems (\ref{auxiliaryproblem}).

\bl \lab{Lemma 2.8}
Assume $(V_1)-(V_3)$ and $(f_1)-(f_3)$ hold, if $\max\{2,2_s^*-2\}<p<2_s^*$, then
for almost every $\la \in [\frac{1}2,1]$, $I_\la$ has a positive critical point $u_\la$ satisfying $\|u_\la\|\geq \delta$ where $\delta>0$ independent of $\la$.
\el
\bp From Lemma \ref {Lemma 2.2} and Lemma \ref{Lemma 2.3}, there exists a bounded sequence $\{u_n\}\subset H^s(\RN),u_n\geq0$ and $0<c_\la<\frac{s}{N\la^{\frac{N-2s}{2s}}}S_s^{\frac{N}{2s}}$, such that
\begin{align*}
I_\la(u_n)\to c_\la,I'_\la(u_n)\to 0.
\end{align*}
Then $u_n\to u_\la\geq 0$ weakly in $H^s(\RN)$. It is obvious that $u_\la$ is a critical point of $I_\la$.

\noindent Now, we claim $u_\la\neq 0$. If $u_\la=0$, from Proposition \ref{proposition1}, we can deduce that  $k>0$ since $c>0$, and \be\lab{clambda1}
c_\la=\sum_{j=1}^k I_\la^\iy(v_\la^j)\geq m_\la^\iy:=\inf\{I_\la^\iy(u):u\in H^s(\RN),u\neq 0,I_\la^{\iy'}(u)=0\},
\ee
where $I_\la^{\iy'}(v_\la^j)=0(j=1,2,\cdots,k)$. On the other hand, we infer that
\be\lab{clamba2}
c_\la<m_\la^\iy,
\ee
which is contradictory to (\ref{clambda1}) and then the claim is true.

From Lemma \ref{Lemma 2.6}, let $v_\la$  be the least energy solution of
\begin{align*}
(-\DD)^su+V(\iy)u=\la f(u).
\end{align*}
By Lemma \ref{Lemma 2.5}, there exists $\g_\la(t)$ satisfying  $\g_\la(0)=0,I_\la^\iy(\g_\la(1))<0,v_\la\in \g_\la[0,1]$ and
\begin{align*}
\max_{t\in[0,1]}I_\la^\iy(\g_\la(t))=I_\la^\iy(v_\la)= m_\la^\iy.
\end{align*}
By $(V_2)$, we have
\begin{align*}
I_\la(\g_\la(t))<I_\la^\iy(\g_\la(t)), \forall t\in [0,1],
\end{align*}
and it follows from the definition of $c_\la$ that
\begin{align*}
c_\la\leq \max_{t\in[0,1]}I_\la(\g_\la(t))<\max_{t\in[0,1]}I_\la^\iy(\g_\la(t))= m_\la^\iy.
\end{align*}

If $(f_1)-(f_2)$ hold, by the same argument as that in Lemma \ref{energyest}, there exists a constant $\delta_0 >0$ independent of $\la$ such that $\int_{\RN}|(-\DD)^{\frac{s}2}u_\la|^2\geq \delta_0$ since $u_\la\neq 0$. Thus, there exists a $\delta>0$ independent of $\la$ such that $\|u_\la\|\geq \delta$. The proof is finished.
\ep

\s{The proof of Theorem \ref{Theorem 1.1}}
\renewcommand{\theequation}{4.\arabic{equation}}
Lemma \ref{Lemma 2.8} shows that for almost every $\la\in [\frac{1}2,1]$, $I_{\la}(u)$ has a positive critical point $u_{\la}$. Thus we get a critical point sequence $\{u_{\la}\}$ satisfying $I'_{\la}(u_{\la})=0$. In the following, we first show that $\{u_\la\}$ is a BPS sequence of $I$  and then prove the convergence of $\{u_\la\}$ as $\la\to 1$. By analyzing the properties of  minimizing sequence, we complete the proof of the existence of ground state solutions of (\ref{p1}).

First, we show that the uniform boundedness of $\{u_\la\}$.
\bo \lab{proposition2} Assume $(V_1)-(V_3)$ and $(f_1)-(f_3)$ hold, if $\max\{2,2_s^*-2\}<p<2_s^*$, then $\{u_{\la}\}$ is bounded uniformly and there exists $\delta>0$ independent of $\la$ such that $I_\la(u_\la)\geq \delta$.
\eo
\bp
Since $u_\la$ is the critical point of $I_\la(u)$, from the Pohoz\v{a}ev identity (\ref{pohozaev identity}), we have
\be
\label{ilambdaest}
I_\la(u_\la)=\frac{s}N\int_{\RN}|(-\DD)^{\frac{s}2}u_\la|^2-\frac{1}{2N}\int_{\RN}\lan\na V(x),x\ran |u_\la|^2.
\ee
From Proposition \ref{proposition1}, $I_\la(u_\la)\leq c_\la\leq c_{\frac{1}2}$ for any $\la\in [\frac{1}2,1]$. By H\"{o}lder inequality and Sobolev embedding theorem,
\begin{align*}
\int_{\RN}|(-\DD)^{\frac{s}2}u_\la|^2&=\frac{N}s I_\la(u_\la)+\frac{1}{2s}\int_{\RN}\lan\na V(x),x\ran |u_\la|^2\\
&\leq \frac{N}s c_{\frac{1}2}+\frac{1}{2sS_s}\|\max\{\lan\na V(x),x\ran,0\}\|_{L^{\frac{N}{2s}}}\int_{\RN}|(-\DD)^{\frac{s}2}u_\la|^2.
\end{align*}
$(V_3)$ implies that $\int_{\RN}|(-\DD)^{\frac{s}2}u_\la|^2$ is bounded uniformly independent of $\la$. Next, we show that $\|u_\la\|_{L^2}$ is bounded uniformly independent of $\la$.
From $I'_\la(u_\la)u_\la=0$, we have $\int_{\RN}|(-\DD)^{\frac{s}2}u_\la|^2+\int_{\RN}V(x)|u_\la|^2=\la\int_{\RN}f(u_\la)u_\la$. Then, by $(f_1)-(f_2)$,
\begin{align*}
V_0\int_{\RN}|u_\la|^2&\leq\int_{\RN}|(-\DD)^{\frac{s}2}u_\la|^2+\int_{\RN}V(x)|u_\la|^2\\
&\leq \la\e\int_{\RN}|u_\la|^2+\la C(\e)\int_{\RN}|u_\la|^{2_s^*}\\
&\leq \e\int_{\RN}|u_\la|^2+C(\e) \left|\int_{\RN}|(-\DD)^{\frac{s}2}u_\la|^2\right|^{\frac{2_s^*}2}.
\end{align*}
Therefore, $\|u_\la\|_{L^2}$  is bounded uniformly. Now, we prove that $I_\la(u_\la)\geq \delta>0$. From Lemma \ref{Lemma 2.8}, there exists $\delta_0>0$ independent of $\la$ such that $\|u_\la\|\geq \delta_0$. On the other hand,
\begin{align*}
I_\la(u_\la)&\geq \frac{s}N\int_{\RN}|(-\DD)^{\frac{s}2}u_\la|^2-\frac{1}{2N}\int_{\RN}\max\{\lan\na V(x),x\ran,0\} |u_\la|^2\\
&\geq \frac{s}N\int_{\RN}|(-\DD)^{\frac{s}2}u_\la|^2-\frac{1}{2NS_s}\|\max\{\lan\na V(x),x\ran,0\}\|_{L^{\frac{N}{2s}}}\int_{\RN}|(-\DD)^{\frac{s}2}u_\la|^2.
\end{align*}
$(V_3)$ implies that there exists $\delta>0$ independent of $\la$ such that
\be\lab{ilambda1}
I_\la(u_\la)\geq \delta
\ee
The proof is finished.
\ep

In the following, we denote $u_\la$ by $u_{\la_j}$ and let $\la_j\to 1$ as $j\to \iy$.
\bl\lab{Lemma 3.1}
Assume $(V_1)-(V_3)$ and $(f_1)-(f_3)$ hold, if $\max\{2,2_s^*-2\}<p<2_s^*$, then the sequence $\{u_{\la_j}\}$ is a BPS sequence for $I$ satisfying $\limsup_{j\to \iy}I(u_{\la_j})\leq c_1$ and $\|u_{\la_j}\|\nrightarrow 0$.
\el
\bp
From Lemma \ref{Lemma 2.8}, we have $\|u_{\la_j}\|\nrightarrow 0$.
 It follows from Proposition \ref{proposition2} that $\|u_{\la_j}\|$ is bounded uniformly, and consequently $\int_{\RN} F(u_{\la_j})$ is bounded by $(f_1)-(f_2)$.  The property $(iii)$ in Proposition \ref{proposition1} shows that $I_{\la_j}(u_{\la_j})\leq c_{\la_j}$ for any $u_{\la_j}$.
Thus, together with
\be \lab{iest}I(u_{\la_j})=I_{\la_j}(u_{\la_j})+(\la_j-1)\int_{\RN} F(u_{\la_j}),\ee
we obtain $\limsup_{j\to \iy}I(u_{\la_j})\leq c_1$ and $I'(u_{\la_j})\to 0$.
\ep

{\bf The completion of the proof of the Theorem \ref{Theorem 1.1}.}
\bp
From Lemma \ref{Lemma 3.1}, together with (\ref{ilambda1}) and (\ref{iest}), we get that there exists a subsequence still denoted by $\{u_{\la_j}\}$ satisfying
\begin{align*}
(i)\ \ \{u_{\la_j}\} \ \ is\ \ bounded,\ \ (ii)\ \ I(u_{\la_j})\to c\leq c_1,\ \ (iii)\ \ I'(u_{\la_j})\to 0,
\end{align*}
where $c>0$. That is to say, there exists a BPS sequence $\{u_{\la_j}\}$ satisfying the assumptions of Lemma \ref{Lemma 2.8} for $\la=1$. Thus, there exists a nontrivial critical point $u_0$ for $I$ satisfying $I(u_0)\leq c_1$.

In the following, we show the existence of a ground state solution. Let
\begin{align*}
m=\inf\{I(u):u\in H^s(\RN),u\neq 0,I'(u)=0\}.
\end{align*}
Obviously, $m\leq I(u_0)\leq c_1=\frac{s}NS_s^{\frac{N}{2s}}$. Set $\{u_n\}$ be a sequence of nontrivial critical points of $I$ satisfying $I(u_n)\to m$.  Since $I(u_n)$ is bounded, similar proof as that in Proposition \ref{proposition2} for $\la=1$, we get that $\{u_n\}$ is bounded uniformly and there exists $\delta>0$ such that $I(u_n)\geq \delta>0$. Thus $m>0$. So, $\{u_n\}$ is a BPS sequence
satisfying the following conditions,
\begin{align*}
(i)\ \ \{u_n\} \ \ is\ \ bounded,\ \ (ii)\ \ I(u_n)\to m\leq c_1,\ \ (iii)\ \ I'(u_n)=0,
\end{align*}
 From Proposition \ref{proposition1}, there exists $\ti u$ such that $I'(\ti u)=0$ and $I(\ti u)\leq m$.

Now, we claim $\ti u\neq 0$.
Otherwise, $\ti u=0$. Then, by Proposition \ref{proposition1}, we have
\begin{align*}
m=\sum_{j=1}^k I^\iy(w^j)\geq m^\iy:=\inf\{I^\iy(u):u\in H^s(\RN),u\neq 0,I^{\iy'}(u)=0\}
\end{align*}
for $k>0$ and $w^j(j=1,2,\cdots,k)$ are the critical points of $I^\iy$.
On the other hand, similar argument as that in Lemma \ref{Lemma 2.8}, there exists  $\g(t)$ such that
\begin{align*}
\max_{t\in[0,1]}I^\iy(\g(t))= m^\iy.
\end{align*}
From the definition of $c_1$, we obtain $m\leq c_1 \leq \max_{t\in[0,1]}I(\g(t))$. By $(V_2)$, we get
\begin{align*}
m\leq c_1< m^\iy,
\end{align*}
which is a contradiction. Thus, the claim is true.
Then $I(\ti u)\geq m$ since $I'(\ti u)=0$ and $\ti u\neq 0$. So, there exists a critical point $\ti u\neq 0$ such that $I(\ti u)=m$. The proof is completed.
\ep

\s*{Acknowledgements}
This work is supported by the National Natural Science Foundation of China (11271364).


\begin{thebibliography}{9999}


\bibitem{G.Alberti} G. Alberti, G. Bouchitt\'{e}, P. Seppecher, Phase transition with the line-tenstion effect. Arch. Ration. Mech. Anal. 144, 1-46(1998).

\bibitem{Alvesfrac} C. O. Alves, O. H. Miyagaki, Existence and concentration of solution for a class of
fractional elliptic equation in $\mathbb{R}^N$ via penalization method, Calc. Var. Partial Differential Equations. 55:47(2016).

\bibitem{B.Barriosa}B. Barriosa, E. Coloradoc, R. Servadeid, F. Soriaa, A critical fractional equation with concave-convex power nonlinearities, Ann. I. H. Poincar\'{e}-AN. 32, 875-900(2015).


\bibitem{Barriosa1} B. Barriosa, E. Colorado, A. de Pablo, U. S¨¢nchez, On some critical problems for the fractional Laplacian operator, J. Differential Equations. 252, 6133-6162 (2012).


\bibitem{Berestycki} H. Berestycki and P. Lions, Nonlinear scalar field equations I. Existence of a ground state, Arch. Ration. Mech. Anal. 82, 90-117(1990).

\bibitem{bisci} G. M. Bisci, V. Radulescu, Ground state solutions of scalar field fractional Schr\"{o}dinger equation, Calc. Var. Partial Differential Equatioans. 54, 2985-3008(2015).

\bibitem{Cabre} X. Cabr¨¦, J. Tan, Positive solutions of nonlinear problems involving the square root of the Laplacian, Adv. Math. 224,2052-2093(2010).

\bibitem{L.A.Caffarelli} L. A. Caffarelli, S. Salsa, L. Silvestre, Regularity estimates for the solution and the free boundary of the obstacle problem for the fractional Laplacian, Invent. Math.  171(2),425-461(2008).

\bibitem{Caffarelli} L. Caffarelli, L. Silvestre, An extension problem related to the fractional Laplacian, Commun. PDEs. 32, 1245-1260(2007).

\bibitem{changxiaojun} X. J. Chang, Z. Q. Wang, Nodal and multiple solutions of nonlinear problems involving the fractiona Laplacian. J. Differential Equations. 256,2965-2992(2014).

\bibitem{changxiaojun2} X. J. Chang, Z. Q. Wang, Ground state of scalar field equations involving a fractional Laplacian with general nonlinearity, Nonlinearity.  26,479-494(2013).

\bibitem{chenzheng}  G. Chen, Y. Zheng, Concentration phenomena for fractional nonlinear Schr\"{o}dinger equations, Commun. Pure Appl. Anal. 13, 2359-2376(2014).

\bibitem{A.Cotsiolis} A. Cotsiolis and N. K. Tavoularis, Best constants for Sobolev inequalities for higher order fractional derivatives, J. Math. Anal. Appl. 295, 225-236(2004).

\bibitem{Davila} J. D\'{a}vila, M. del Pino, J. Wei, Concentrating standing waves for the fractional nonlinear Schr\"{o}ding equation, J. Differ. Equ. 256(2), 858-892(2014).

\bibitem{FMM}  M. M. Fall,  F. Mahmoudi, E. Valdinoci, Ground states and concentration phenomena for the fractional Schr\"{o}dinger equations, Nonlinearity. 28, 1937-1961(2015).


\bibitem{Q.Y.Guan} Q. Y. Guan, Z. M. Ma, Boundary problems for fractional Laplacians, Stoch. Dyn. 593,385-424(2005).

\bibitem{xiaominghe} X. M. He, W. M. Zou, Existence and concentration result for the fractional Schr\"{o}dinger equations with critical nonlinearities, Calc. Var. Partial Differential Equations. 55(4):91(2016).

\bibitem{L.Jeanjean} L. Jeanjean, On the existence of bounded Palais-Smale sequence and application to a Landesman-Lazer-type problem set on $\RN$, Proc. Roy. Soc. Edinburgh. 129A, 787-809(1999).

\bibitem{L.JeanjeanandKazunagaTanaka}L. Jeanjean and K. Tanaka, A Positive solution for a nonlinear Schr\"{o}dinger equation on $\RN$. Indiana. Univ. Math. J. 54, 443-464(2005).

\bibitem{JeanTa} L. Jeanjean and K. Tanaka, A remark on least energy solutions in $\mathbb {R^N}$, Proc. Amer. Math. Soc. 13, 2399-2408(2002).

\bibitem{NLaskin1} N. Laskin, Fractional Schr\"{o}dinger equation, Phys. Rev. E66, 05618(2002).

\bibitem{NLaskin2} N. Laskin, Fractional quantum mechanics, Phys. Rev. E62, 3135(2000).

\bibitem{Laskin} N. Laskin, Fractional quantum mechanics ans L\'{e}vy path integrals, Phys. Lett. A 268(4-6), 298-305(2000).

\bibitem{lionsembedding} P. L. Lions, Sym\'{e}trie et compacit\'{e} dans les espaces de Sobolev, J. Funct. Analysis. 49, 315-334(1982).


\bibitem{Nezza} E. Di Nezza, G. Palatucci and E. Valdinoci, Hitchiker's guide to the fractional Sobolev space. Bull. Sci. Math. 136,  521-573(2012)

\bibitem{lionslemma} S. Secchi, Ground state solutions for nonlinear fractional Schr\"{o}ding equations in $\RN$. J. Math. Phys. 54,  031501(2013).

\bibitem{seok} J. Seok, Spike-layer solutions to nonlinear fractional Schr\"{o}ding equations with almost optimal nonlinearities,
Electron. J. Differential Equations. 196, 1-19(2015).

\bibitem{shangxzhangj} X. D. Shang, J. H. Zhang, Concentrating solutions of nonlinear fractional Schr\"{o}dinger equation with potentials, J. Differ. Equ. 258(4), 1106-1128(2015).

\bibitem{Silvestre}L. Silvestre, Regularity of the obstacle problem for a fractional power of the Laplace ooperator, Comm. Pure Apple. Math. 60(1), 67-112(2007).

\bibitem{Y.Sire} Y. Sire, E. Valdinoci, Fractional Laplacian phase transtion and boundary reactions: a geometric inequality and a symmetry result. J. Funct. Anal. 256(6), 1842-1864(2009).

\bibitem{simonesecchi} S. Secchi, On fractional Schr\"{o}dinger equations in $\RN$ without the Ambrosetti-Rabinowitz condition. Topological Methods in Nonlinear Analysis. 47(1), 19-41(2016).

\bibitem{J.Tan} J. Tan, The Brezis-Nirenberg type problem involving the square root of the Laplacian, Calc. Var. 42, 21-41 (2011).

\bibitem{tengkaiming} K. M. Teng, X. M. He, Ground state solutions for fractional Schr\"{o}dinger equations with critical Sobolev exponent, Commun. Pure Apple. Anal. 16, 991-1008(2016).

\bibitem{weiyuanhong} Y. H. wei, X. F. Su, Multiplicity of solutions for non-local elliptic equations driven by the fractional Laplacian, Calc. Var. 52, 95-124(2015).

\bibitem{J.Zhang} J. Zhang, W. M. Zou, The critical case for a Berestycki-Lions theorem, Science China Mathematics. 57(3), 541-554(2014).

\bibitem{brezis-lieb} J. J. Zhang, J. M. do O and M. Squassina, Schr\"{o}dinger-Poisson systems with a general critical nonlinearity. Commun. Contemp. Math. (published online 2016).

\bibitem{zjjjm} J. J. Zhang, J. M. do O and M. Squassina, Fractional Schr\"{o}dinger-Poisson systems with a general subcritical or critical nonlinearity, Adv. Nonlinear Stud. 16(1), 15-30(2016).

\bibitem{zjjzwm2} J. J. Zhang and W. M. Zou, A Berestycki-Lions Theorem revisited, Commun. Contemp. Math. 14(5), 1250033(2012).

\bibitem{zjjzwm} J. J. Zhang, W.M. Zou, Solutions concentrating around the saddle points of the potential for critical schr\"{o}dinger equations, Calc. Var. Partial Differential Equations. 54(4), 4119-4142(2015).

\bibitem{zhangxia} X. Zhang, B. L. Zhang, D. repov\u{s}, Existence and symmetry of solutions for critical fractional Schr\"{o}dinger equations with bounded potentials. Nonlinear Analysis. 142, 48-68(2016).


\end{thebibliography}
\end{document}